\input amstex
\documentstyle{amsppt}
\magnification=\magstep1
 \hsize 13cm \vsize 18.35cm \pageno=1
\loadbold \loadmsam
    \loadmsbm
    \UseAMSsymbols
\topmatter
\NoRunningHeads
\title A note on   the $q$-Genocchi  numbers and polynomials
\endtitle
\author
  Taekyun Kim
\endauthor
 \keywords $p$-adic $q$-integrals, Genocchi numbers, q-Genocchi numbers,
 polynomials, sums of powers
\endkeywords

\abstract
 In this paper we discuss the new concept of the $q$-extension of Genocchi
 numbers and give the some relations between $q$-Genocchi
 polynomials and $q$-Euler numbers.
\endabstract
%\thanks  2000 AMS Subject Classification:
%\newline keywords and phrases :
%\newline  This paper is supported by  Jangjeon Mathematical
%Society
%\endthanks
\endtopmatter

\document

{\bf\centerline {\S 1. Introduction}}

 \vskip 20pt
A number given by the generating function
$$\frac{2t}{e^t +1}=\sum_{n=0}^{\infty} G_n \frac{t^n}{n!}. \tag 1$$
It satisfies  $G_1=1, G_3=G_4=G_7=\cdots =0, $ and even
coefficients are given by $G_{2m}=2(1-2^{2n})B_{2n}=2nE_{2n-1}(0)
,$ where $B_n$ is a Bernoulli number and $E_n(x)$ is an Euler
polynomial. The first few Genocchi numbers for $n=2,4,\cdots$ are
$-1, -3, 17, -155, 2073, \cdots .$ The Euler polynomials are well
known as
$$\frac{2}{e^t +1}e^{xt}=\sum_{n=0}^{\infty}E_n(x)\frac{t^n}{n!},
\text{ see [1, 3, 7, 8, 9]. } \tag 2$$ By (1) and (2) we easily
see that
$$E_n(x)=\sum_{k=0}^n \binom{n}{k}\frac{G_{k+1}}{k+1} x^{n-k},
\text{ where $\binom{n}{k}=\frac{n (n-1)\cdots(n-k+1)}{k!}$, cf.
[4, 5, 6]. }\tag 3$$ For $m,n \geq 1$ and $m$ odd, we have
$$(n^m -n)G_m =\sum_{k=1}^{m-1}\binom{m}{k}n^kG_k Z_{m-k}(n-1),
\tag 4$$ where $Z_m(n)=1^m-2^m+3^m-\cdots+(-1)^{n-1}n^m$, see[3].
From (10 we derive
$$2t =\sum_{n=0}^{\infty}\left((G+1)^n+G_n \right) \frac{t^n}{n!},$$
where we use the technique method notation by replacing $G^m $ by
$G_m (m \geq 0),$ symbolically. By (5), we see that
$$
G_0 =0, \quad (G +1)^n -G_n = \cases
2 & \text{if\   $n=1$}\\
0 & \text{if \  $n>1$.}
\endcases\tag 6$$
Let $p$ be a fixed odd prime, and let $\Bbb C_p$ denote the
$p$-adic completion of the algebraic closure of $\Bbb Q_p
$(=$p$-adic number field). For $d$ a fixed positive integer with
$(p,d)=1$, let
$$\split
& X=X_d = \lim_{\overleftarrow{N} } \Bbb Z/ dp^N \Bbb Z ,\cr & \
X_1 = \Bbb Z_p , \cr  & X^\ast = \underset {{0<a<d p}\atop
{(a,p)=1}}\to {\cup} (a+ dp \Bbb Z_p ), \cr & a+d p^N \Bbb Z_p =\{
x\in X | x \equiv a \pmod{dp^n}\},\endsplit$$ where $a\in \Bbb Z$
lies in $0\leq a < d p^N$.

Ordinary $q$-calculus is now very well understood from many
different point of views. Let us consider a complex number $q\in
\Bbb C$ with $|q|<1$ (or $q\in \Bbb C_p$ with $|1-q|_p <
p^{-\frac{1}{p-1}}$) as an indeterminate. The $q$-basic numbers
are defined by
$$
[x]_q  = \dfrac{q^x -1}{q -1} = 1+ q + q^2 +\cdots + q^{x-1},
$$
and

$$
[x]_{-q}  = \dfrac{-(-q)^x +1}{q +1} = 1- q + q^2 +\cdots +
(-q)^{x-1}.
$$
We say that $f$ is a uniformly differentiable function at a point
$a \in\Bbb Z_p $ and denote this property by $f\in UD(\Bbb Z_p )$,
if the difference quotients
$$
F_f (x,y) = \dfrac{f(x) -f(y)}{x-y}
$$
have a limit $l=f^\prime (a)$ as $(x,y) \to (a,a)$.

For $f\in UD(\Bbb Z_p )$, let us start with the expression

$$\eqalignno{ & \dfrac{1}{[p^N ]_q} \sum_{0\leq j < p^N} q^j f(j) =\sum_{0\leq j < p^N} f(j)
\mu_q (j +p^N \Bbb Z_p ), }
$$
representing a $q$-analogue of Riemann sums for $f$, cf. [5]. The
integral of $f$ on $\Bbb Z_p$ will be defined as limit ($n \to
\infty$) of those sums, when it exists. The $p$-adic $q$-integral
of the function $f\in UD(\Bbb Z_p )$ is defined by
$$ \eqalignno{ &
I_q (f) =\int_{\Bbb Z_p }f(x) d\mu_q (x) = \lim_{N\to \infty}
\dfrac{1}{[p^N ]_q} \sum_{0\leq x < p^N} f(x) q^x ,\quad
\text{(see [5, 10, 11, 12] )}. }$$ In the previous paper [4]
author constructed the $q$-extension of Euler polynomials in the
fermionic sense of $p$-adic $q$-integral at $q=-1$ as follows:
$$E_{n,q}(x)=\int_{\Bbb Z_p}[t+x]_q^n d\mu_{-q}(t), \text{ where $\mu_{-q}(x+p^N\Bbb Z_p)=
\lim_{q\rightarrow -q}\mu_{q}(x+p^N \Bbb Z_p) .$}\tag 7$$ From
(7), we note that
$$E_{n,q}(x)=\frac{[2]_q}{(1-q)^n}\sum_{l=0}^n
\binom{n}{l}\frac{(-1)^l}{1+q^{l+1}}q^{lx}, \text{ see [4].}\tag
8$$ The $q$-extension of Genocchi numbers is defined as
$$g_q^{*}(t)=[2]_qt\sum_{n=0}^{\infty}(-1)^nq^n
e^{[n]_qt}=\sum_{n=0}^{\infty}G_{n,q}^{*}\frac{t^n}{n!}, \text{
see [4]}. \tag9$$ The following formula is well known in [4, 7]:
$$E_{n,q}(x)=\sum_{k=0}^n\binom{n}{k}[x]_q^{n-k}q^{kx}\frac{G_{k+1,q}^{*}}{k+1}.\tag 10$$
The modified $q$-Euler numbers are defined  as
$$\eqalignno{ &
\xi_{0,q} =\frac{[2]_q}{2},\quad  (q\xi +1)^k +\xi_{k,q} = \cases
[2]_q & \text{if\   $k=0$}\\
0 & \text{if \  $k\neq0$,}
\endcases }$$ with the usual convention of replacing $\xi^i$ by $\xi_{i,q}$, see
[10]. Thus, we derive the generating function of $\xi_{n,q}$ as
follows:
$$F_q(t)=[2]_q\sum_{k=0}^{\infty}(-1)^ke^{[k]_qt}=\sum_{n=0}^{\infty}\xi_{n,q}\frac{t^n}{n!}.\tag11$$
Now we also consider the $q$-Euler polynomials $\xi_{n,q}(x)$ as
$$F_q(t,x)=[2]_q\sum_{k=0}^{\infty}(-1)^ke^{[k+x]_qt}=\sum_{n=0}^{\infty}\xi_{n,q}(x)\frac{t^n}{n!}.\tag12$$
From (12) we note that
$$\xi_{n,q}(x)=\sum_{l=0}^n
\binom{n}{l}\xi_{l,q}q^{lx}[x]_q^{n-l}, \text{ see [10]}. \tag13$$
In the recent several authors studied the $q$-extension of
Genocchi numbers and polynomials (see[ 1, 2, 5, 6, 7, 12]). In
this paper we discuss the new concept of the $q$-extension of
Genocchi numbers and give the some relations between $q$-Genocchi
numbers and $q$-Euler numbers.

\vskip 20pt

{\bf\centerline {\S 2. $q$-extension of Genocchi numbers}}

 \vskip 20pt

In this section we assume that $q\in \Bbb C$ with $|q|<1.$ Now we
consider the $q$-extension of Genocchi numbers as follows:
$$g_q(t)=[2]_qt\sum_{k=0}^{\infty}(-1)^ke^{[k]_qt}=\sum_{n=0}^{\infty}G_{n,q}\frac{t^n}{n!}.\tag14$$
In Eq.(14), it is easy to show that $\lim_{q\rightarrow 1}
g_q(t)=\frac{2t}{e^t+1}=\sum_{n=0}^{\infty}G_n \frac{t^n}{n!}.$
From (14) we derive

$$\eqalignno{g_q(t)&=[2]_q
t\sum_{k=0}^{\infty}(-1)^k\sum_{m=0}^{\infty}[k]_q^m\frac{t^m}{m!}=[2]_q\sum_{k=0}^{\infty}
(-1)6k\sum_{m=1}^{\infty}m[k]_q^{m-1}\frac{t^m}{m!} \cr&=
[2]_q\sum_{k=0}^{\infty}(-1)^k\sum_{m=0}^{\infty}m[k]_q^{m-1}\frac{t^m}{m!}.
&(15)}
$$
By (15), we easily see that
$$g_q(t)=[2]_q\sum_{m=0}^{\infty}\left(m
\left(\frac{1}{1-q}\right)^{m-1}\sum_{l=0}^{m-1}\binom{m-1}{l}(-1)^l\frac{1}{1+q^l}\right)\frac{t^m}{m!}.\tag16$$
From (14) and (16) we note that
$$\sum_{m=0}^{\infty}G_{m,q}\frac{t^m}{m!}=\sum_{m=0}^{\infty}\left(m[2]_q\left(\frac{1}{1-q}\right)^{m-1}\sum_{l=0}^{m-1}
\binom{m-1}{l}\frac{(-1)^l}{1+q^l}\right)\frac{t^m}{m!}.\tag17$$
By comparing the coefficients on both sides in Eq.(17), we have
the following theorem:

\proclaim{ Theorem 1} For $ m\geq 0$ we have
$$G_{m,q}=m[2]_q\left(\frac{1}{1-q}\right)^{m-1}
\sum_{l=0}^{m-1}\binom{m-1}{l}\frac{(-1)^l}{1+q^l}. $$
\endproclaim
From Theorem 1, we easily derive the following corollary:
\proclaim{ Corollary 2} For $n \in \Bbb N$ we have
$$\eqalignno{ &
G_{0,q} =0,\quad  (qG+1)^k +G_{k,q} = \cases\frac{[2]_q^2}{2}
 & \text{if\   $k=1$},\\
0 & \text{if \  $k>1$,}
\endcases }$$ with the usual convention of replacing $G^{i}$ by
$G_{i,q}.$
\endproclaim
Remark. We note that Corollary 2 is the $q$-extension of Eq.(6).
By (10) and Corollary 2, we obtain the following theorem:
\proclaim{ Theorem 3} For $n\in \Bbb N$ we have
$$\xi_{n,q}=\frac{G_{n+1,q}}{n+1}.$$
\endproclaim
From (12) we derive

$$\eqalignno{F_q(x,t)&=[2]_q\sum_{n=0}^{\infty}(-1)^ne^{[n+x]_qt}=q^xt\frac{[2]_q}{q^x
t}e^{[x]_qt}\sum_{n=0}^{\infty}(-1)^ne^{q^x [n]_qt} \cr & =
e^{[x]_qt}\sum_{n=0}^{\infty}q^{nx}\frac{G_{n+1,q}}{n+1}\frac{t^n}{n!}
=\sum_{n=0}^{\infty}\left(\sum_{k=0}^n\binom{n}{k}[x]_q^{n-k}q^{kx}\frac{G_{k+1,q}}{k+1}\right)\frac{t^n}{n!}.
&(18) }$$ By (18), we easily see that
$$\xi_{n,q}(x)=\sum_{k=0}^n\binom{n}{k}[x]_q^{n-k}q^{kx}\frac{G_{k+1,q}}{k+1}.$$
This formula can be considered as the $q$-extension of Eq.(3). Let
us consider the $q$-analogue of Genocchi polynomials as follows:
$$g_q(x,t)=[2]_q
t\sum_{k=0}^{\infty}(-1)^ke^{[k+x]_qt}=\sum_{n=0}^{\infty}G_{n,q}(x)\frac{t^n}{n!}.
\tag19$$ Thus, we note that $\lim_{q\rightarrow
1}g_q(x,t)=\frac{2t}{e^t
+1}e^{xt}=\sum_{n=0}^{\infty}G_n(x)\frac{t^n}{n!}.$ From Eq.(19),
we easily derive
$$
G_{n,q}(x)=[2]_qn\left(\frac{1}{1-q}\right)^{n-1}\sum_{l=0}^{n-1}\frac{(-1)^l}{1+q^l}q^{lx}\binom{n-1}{l}.\tag20$$

By (19) we also see that

$$\eqalignno{
\sum_{n=0}^{\infty}G_{n,q}(x)\frac{t^n}{n!} &=[2]_q
t\sum_{k=0}^{\infty}(-1)^ke^{[k+x]_qt}=[2]_qt\sum_{a=0}^{m-1}(-1)^a\sum_{k=0}^{\infty}(-1)^ke^{[k+\frac{a+x}{m}]_{q^m}
[m]_qt} \cr
&=\frac{[2]_q}{[m]_q[2]_{q^m}}\sum_{a=0}^{m-1}(-1)^a\left([m]_qt[2]_{q^m}\sum_{k=0}^{\infty}(-1)^k
e^{[m]_qt[k+\frac{a+x}{m}]_{q^m}}\right)\cr
&=\sum_{n=0}^{\infty}\left(\frac{[2]_q}{[m]_q[2]_{q^m}}\sum_{a=0}^{m-1}(-1)^a
[m]_q^n G_{n,q^m}\left(\frac{x+a}{m}\right)
\right)\frac{t^n}{n!}\cr&=\sum_{n=0}^{\infty}\left(\frac{[2]_q}{[2]_{q^m}}[m]_q^{n-1}\sum_{a=0}^{m-1}(-1)^a
G_{n,q^{m}}\left(\frac{x+a}{m}\right)\right)\frac{t^n}{n!}, \text{
where $m\in\Bbb N$ odd.}  }$$ Therefore we obtain the following
theorem: \proclaim{ Theorem 4} Let $m(=odd)\in\Bbb N$. Then we
have the distribution of the $q$-Genocchi polynomials as follows:
$$G_{n,q}(x)=\frac{[2]_q}{[2]_{q^m}}[m]_q^{n-1}\sum_{a=0}^{m-1}(-1)^aG_{n,q^m}\left(\frac{x+a}{m}\right),$$
where $n$ is positive integer.
\endproclaim
Theorem 4 will be used to construct the $p$-adic $q$-Genocchi
measures which will be treated in the next section. Let $\chi$ be
a primitive Dirichlet character with a conductor $d(=odd)\in \Bbb
N$. Then the generalized $q$-Genocchi numbers attached to $\chi$
are defined as
$$g_{\chi,q}(t)=[2]_q t\sum_{a=0}^{d-1}\chi(n)(-1)^ne^{[n]_q
t}=\sum_{n=0}^{\infty}G_{n,\chi,q}\frac{t^n}{n!}. \tag21$$ From
(21) we derive
$$G_{n,\chi,q}=\frac{[2]_q}{[2]_{q^d}}[d]_q^{n-1}\sum_{a=0}^{d-1}(-1)^a\chi(a)G_{n,q^d}(\frac{a}{d}).\tag22$$

{\bf\centerline {\S 3. $p$-adic $q$-Genocchi measures}}

 \vskip 20pt

In this section we assume that $q\in\Bbb C_p$ with
$|1-q|_p<p^{-\frac{1}{p-1}}$ so that $q^x= \exp (x \log q) $. Let
$\chi$ be a primitive Dirichlet's character with a conductor
$d(=odd)\in\Bbb N.$ For any positive integers $N, k$ and
$d(=odd)$, let $\mu_k=\mu_{k,q; G}$ be defined as
$$\mu_k(a +dp^N\Bbb
Z_p)=(-1)^a[dp^N]_q^{k-1}\frac{[2]_q}{[2]_{q^{dp^N}}}G_{k,q^{dp^N}}(\frac{a}{dp^N}).
\tag23$$ By using Theorem 4 and (23), we show that
$$\sum_{i=0}^{p-1}\mu_k(a+idp^N+dp^{N+1}\Bbb Z_p)=\mu_k(a+dp^N\Bbb
Z_p).$$ Therefore we obtain the following theorem: \proclaim{
Theorem 5} Let $d$ be odd positive integer. For any positive
integers $N, k,$ and let $\mu_k=\mu_{k,q;G}$ be defined as
$$\mu_k(a +dp^N\Bbb
Z_p)=(-1)^a[dp^N]_q^{k-1}\frac{[2]_q}{[2]_{q^{dp^N}}}G_{k,q^{dp^N}}(\frac{a}{dp^N}).$$
Then $\mu_k$ can be extended to a distribution on $X$.
\endproclaim
From the definition of $\mu_k$ and (22) we note that
$$\int_{X}\chi(x)d\mu_k(x)=G_{k,\chi,q}.$$
By (14) and (16), it is not difficult to show that
$$G_{n,q}(x)=\sum_{k=0}^n\binom{n}{k}[x]_q^{n-k}q^{kx}G_{k,q}.
\tag24$$ From (23) and (24) we derive
$$d\mu_k(a)=\lim_{N\rightarrow \infty}\mu_{k}(a+dp^N \Bbb Z_p)
=k[a]_q^{k-1}d\mu_{-q}(a).\tag25$$

Therefore we obtain the following corollary: \proclaim{Corollary
6} Let $k$ be a posotove integer. Then we have
$$G_{k,\chi,q}=\int_{X}\chi(x)d\mu_k(x)=k\int_{X}\chi(x)[x]_q^{k-1}d\mu_{-q}(x).$$
Moreover, $$G_{k,q}=k\int_{X}[x]_q^{k-1}d\mu_{-q}(x).$$
\endproclaim
Remark. In the recent paper( see [1]) Cenkci-Can-Kurt have studied
$q$-Genocchi numbers and polynomials and $p$-adic $q$-Genocchi
measures. Our $q$-Genocchi numbers and polynomials to treat in
this paper are different their $q$-Genocchi numbers and
polynomials.

\vskip 20pt

\quad Taekyun Kim

\quad EECS, Kyungpook National University, Taegu 702-701, S. Korea

\quad e-mail:\text{ tkim$\@$knu.ac.kr; tkim64$\@$hanmail.net}

\enddocument